\documentclass[12pt]{amsart}
\usepackage{a4wide}
\usepackage[utf8]{inputenc}
\usepackage{amsmath}
\usepackage{amsfonts}
\usepackage{amssymb}
\usepackage{amsthm}
\usepackage{subfigure}
\newtheorem{theo}{Theorem}[section]

\theoremstyle{definition}

\theoremstyle{remark}
\newtheorem{rema}[theo]{Remark}

\newcommand{\nwc}{\newcommand}
\nwc{\eps}{\epsilon}
\nwc{\ep}{\epsilon}
\nwc{\vareps}{\varepsilon}
\nwc{\Oph}{\operatorname{Op}_\hbar}
\nwc{\la}{\langle}
\nwc{\ra}{\rangle}

\nwc{\mf}{\mathbf} 
\nwc{\blds}{\boldsymbol} 
\nwc{\ml}{\mathcal} 

\nwc{\defeq}{\stackrel{\rm{def}}{=}}

\nwc{\cE}{\ml{E}}
\nwc{\cN}{\ml{N}}
\nwc{\cO}{\ml{O}}
\nwc{\cP}{\ml{P}}
\nwc{\cU}{\ml{U}}
\nwc{\cV}{\ml{V}}
\nwc{\cW}{\ml{W}}
\nwc{\tU}{\widetilde{U}}
\nwc{\IN}{\mathbb{N}}
\nwc{\IR}{\mathbb{R}}
\nwc{\IZ}{\mathbb{Z}}
\nwc{\IC}{\mathbb{C}}
\nwc{\IT}{\mathbb{T}}
\nwc{\IS}{\mathbb{S}}
\nwc{\tP}{\widetilde{P}}
\nwc{\tPi}{\widetilde{\Pi}}
\nwc{\tV}{\widetilde{V}}
\nwc{\supp}{\operatorname{supp}}
\nwc{\rest}{\restriction}
\nwc{\todo}[1]{$\clubsuit$ {\tt #1}}

\renewcommand{\Im}{\operatorname{Im}}

\begin{document}

\title[Observability and quantum limits for the schr\"odinger equation on $\IS^d$]{Observability and quantum limits\\ for the schr\"odinger equation on $\IS^d$}

\author[Fabricio Maci\`a]{Fabricio Maci\`a}
\author[Gabriel Rivi\`ere]{Gabriel Rivi\`ere}

\address{Universidad Polit\'ecnica de Madrid. DCAIN, ETSI Navales. Avda. Arco de la Victoria s/n. 28040 Madrid, Spain}
\email{Fabricio.Macia@upm.es}
\thanks{FM takes part into the visiting faculty program of ICMAT and is
partially supported by grants ERC Starting Grant 277778 and
MTM2013-41780-P (MEC)}

\address{Laboratoire Paul Painlev\'e (U.M.R. CNRS 8524), U.F.R. de Math\'ematiques, Universit\'e Lille 1, 59655 Villeneuve d'Ascq Cedex, France}
\email{gabriel.riviere@math.univ-lille1.fr}
\thanks{GR is partially supported by the Agence Nationale de la Recherche through the
Labex CEMPI (ANR-11-LABX-0007-01) and the ANR project GeRaSic (ANR-13-BS01-
0007-01)}

\begin{abstract}
In this note, we describe our recent results on semiclassical measures for the Schr\"odinger evolution on Zoll manifolds. We focus on the 
particular case of eigenmodes of the Schr\"odinger operator on the sphere endowed with its canonical metric. We also 
recall the relation of this problem with the observability question from control theory. In particular, we exhibit examples of open sets and potentials on 
the 2-sphere for which observability fails for the evolution problem while it holds for the stationary one. Finally, we give some new results in the case 
where the Radon transform of the potential identically vanishes.
\end{abstract}

\maketitle

\section{Introduction}

Let $\IS^d$ be the sphere of dimension $d\geq 2$ endowed with its canonical metric and let $V$ be a smooth real valued function on $\IS^d$. 
Our goal here is to understand the behavior of the Schr\"odinger eigenfunctions:
\begin{equation}\label{e:eigenf}
\left(-\frac{1}{2}\Delta+V(x)\right)u(x)=\lambda u(x), \quad ||u||_{L^2(\IS^d)}=1,
\end{equation}
in the high-frequency limit $\lambda \rightarrow \infty$. Such functions can be identified with stationary solutions of the 
following Schr\"odinger equation:
\begin{equation}\label{e:nonsemiclassical-schr}
i\partial_tv(t,x)=\left(-\frac{1}{2}\Delta+V(x)\right)v(t,x),\quad v|_{t=0}=u_0\in L^2(\IS^d).
\end{equation}
Solutions of~\eqref{e:nonsemiclassical-schr} encode the position probability density of a quantum particle confined on the surface of the sphere and propagating under the action of the potential $V$. In this note, 
we revisit and extend some of the the results in \cite{MaRi16}, and reinterpret them from the light of control theory for the Schrödinger equation.
\subsection{Controllability and observability}
Let us briefly recall the basics of controllabilty theory for this equation. Fix $\omega$ an open set in 
$\IS^d$ and some final time $T>0$. The controllability problem for~\eqref{e:nonsemiclassical-schr} is the following. Given $u_0$ and $u_1$ in $L^2(\IS^d)$, is 
it possible to find $f(t,x)$ in $L^2([0,T]\times\IS^d)$ such that the solution $\psi(t,x)$ of
\begin{equation}\label{e:schr-control}
i\partial_t\psi(t,x)+\left(\frac{1}{2}\Delta-V(x)\right)\psi(t,x)=\mathbf{1}_{\omega}(x)f(t,x),\quad \psi|_{t=0}=u_0
\end{equation}
satisfies $\psi|_{t=T}=u_1$? In other words, can you drive any $u_0$ to any $u_1$ in time $T$ through the Schr\"odinger evolution by acting only the set 
$\omega$? If this is possible, we say that \emph{the Schr\"odinger equation is controllable in time $T$ on the open set $\omega$.} 

It turns out that the controllability property is equivalent to a stability-type estimate for the solutions to the homogeneous Schrödinger equation~\eqref{e:nonsemiclassical-schr}. The \emph{Schr\"odinger equation is said to be 
observable on the set $\omega$ in time $T>0$} if there exists $C_{\omega, T}>0$ such that
\begin{equation}\label{e:observability}\forall u_0\in L^2(\IS^d),\ \|u_0\|_{L^2(\IS^d)}^2\leq C_{\omega,T}\int_0^T\|v(t,x)\|_{L^2(\omega)}^2dt,\end{equation}
where $v(t,x)$ is the solution to the homogeneous Schrödinger equation~\eqref{e:nonsemiclassical-schr} with initial data $u_0$. It turns out that the 
controllability property for \eqref{e:schr-control} and the observability for \eqref{e:nonsemiclassical-schr} are equivalent notions. The simple proof of this fact is part of the so-called
Hilbert Uniqueness Method~\cite{Li88}. Let us briefly recall it here for the sake of completeness. At the expense of replacing $T$ by $T/2$, the problem reduces to studying the particular case $u_1=0$. One then  considers the operator $\Lambda$ defined by:
$$
\Lambda: L^2((0,T)\times \omega)\ni f \longmapsto \psi_f |_{t=0}\in L^2(\IS^d),
$$ 
where $\psi_f$ is the solution to \eqref{e:schr-control} with control $f$ that satisfies $\psi_f |_{t=T}=0$. The fact that the equation is controllable in time $T$ on the open set $\omega$ is equivalent to the fact that the linear bounded operator $\Lambda$ is onto. This property, in turn, is equivalent to the unique solvability of the adjoint equation with an estimate:
$$\Lambda^*u_0 = f\in\Im \Lambda^*,\quad ||u_0||_{L^2(\IS^d)}^2\leq C ||\Lambda^*u_0||_{L^2((0,T)\times \omega)}^2,$$
by the closed graph theorem. It is straightforward to check that $\Lambda^*u_0=-i\mathbf{1}_\omega v$, where $v$ is the solution to \eqref{e:nonsemiclassical-schr} with initial datum $v_0$ and therefore the result follows with $C_{\omega,T}=C$. 

A remarkable result of Lebeau states that observability (and thus 
control of the Schr\"odinger equation) holds for any 
$T>0$ on the open set $\omega$ provided that the following \emph{geometric control condition} is satisfied~\cite{Le92}:
\begin{equation}\label{e:geometric-control}K_{\omega}:=\left\{\gamma\ \text{closed geodesic of}\ \IS^d: \gamma\cap\omega=\emptyset\right\}=\emptyset.
\end{equation}
Conversely, one can show that, if $K_{\overline{\omega}}\neq\emptyset$, then observability fails \emph{for any choice of $V$ in $\ml{C}^{\infty}(\IS^d;\IR)$} -- see for instance~\cite[Prop.~2.2]{MaRi16}. 
The same result holds if $\IS^d$ is replaced by a Riemannian manifold all whose geodesics are closed (these are called \textit{Zoll manifolds}), see~\cite{MaciaDispersion, MaRi16}. 
Whereas Lebeau's result holds for any compact Riemannian manifold, the Geometric Control Condition is not necessary in general. For instance, observability holds under weaker hypotheses on $\omega$ 
on flat manifolds, see for instance \cite{AFM15, ALM16, AM14, BBZ14, BZ12}, or negatively curved manifolds \cite{AnRiv}. 

\subsection{Observability and Quantum Limits}
When particularized to stationary solutions of \eqref{e:nonsemiclassical-schr}, Lebeau's theorem shows that for every $\omega$ 
satisfying $K_{\omega}=\emptyset$, there exists $C_{\omega}>0$ such that, for every $u$ solution of~\eqref{e:eigenf}, one has
\begin{equation}\label{e:observability-eigenfunction}0<C_{\omega}\leq\int_{\omega}|u(x)|^2\text{vol}(dx),\end{equation}
where $\text{vol}$ is the canonical volume measure on $\IS^d$. Since the constant $C_\omega$ is independent of the frequency, estimate \eqref{e:observability-eigenfunction} provides a restriction on the regions in $\IS^d$ on which the $L^2$-mass of high-frequency eigenfunctions can concentrate. We refer to~\cite{BuZw04,Mi12} for more explicit relations between observability for  eigenmodes (or quasimodes) and for the Schr\"odinger evolution. 

In the case where $V\equiv 0$, eigenfunctions are merely spherical harmonics. Using their explicit expression one can prove that the observability estimate \eqref{e:observability-eigenfunction} fails as soon as $K_{\overline{\omega}}\neq\emptyset$. We refer the reader to the work by Jakobson and 
Zelditch~\cite{JaZe99} for the proof of a stronger result -- see also~\cite{AzM10, Ma08} for 
alternative proofs that extend to other manifolds than the sphere.  

Note that in spite of the fact that the observability estimate for eigenfunctions \eqref{e:observability-eigenfunction} is weaker than the corresponding estimate for time-dependent solutions \eqref{e:observability}, 
the conditions on $\omega$ under which these estimates hold are exactly the same when $V$ vanishes identically on $\IS^d$. In fact, the same phenomenon takes place on the planar disk under a weaker geometric condition: both estimates hold if $\omega$ intersects the boundary on an open set, and fail if $\overline{\omega}$ is strictly contained in the interior of the disk \cite{ALM16}. On the flat torus, both estimates hold for any open set $\omega$, even in the presence of a non-zero potential \cite{AFM15,AM14,BBZ14,BZ12}.

It is therefore natural to ask whether or not estimates \eqref{e:observability} and \eqref{e:observability-eigenfunction} are equivalent, \emph{i.e.} on any compact manifold both estimates hold for the same class of open sets $\omega$. In this note, we answer this question by the negative. Our examples are precisely Schrödinger operators on the sphere with non-constant potentials, or more generally, Laplacians on Zoll manifolds. 

Before stating our results, let us mention that these questions are naturally related to certain problems arising in mathematical 
physics. In fact, consider the set $\mathcal{N}(\infty)$ of probability measures in $\IS^d$ that are obtained as follows. A probability measure $\nu$ 
belongs to $\ml{N}(\infty)$ provided there exists a sequence of eigenfunctions $(u_n)$ :
$$-\frac{1}{2}\Delta u_n+Vu_n=\lambda_n u_n, \quad ||u_n||_{L^2(\IS^d)}=1,$$
with eigenvalues satisfying $\lambda_n\rightarrow\infty$ such that 
$$\lim_{n\rightarrow\infty}\int_{\IS^d} a(x)|u_n|^2(x)\operatorname{vol}(dx)=\int_{\IS^d} a(x)\nu(dx),\quad \text{for every }a\in\ml{C}^0(\IS^d).$$
Measures in $\ml{N}(\infty)$ therefore describe the asymptotic mass distribution sequences of eigenfunctions $(u_n)$ whose corresponding 
eigenvalues tend to infinity. If one integrates these objects against $a=\mathbf{1}_{\omega}$, then one recovers the quantity we were considering before. In quantum mechanics, they describe the probability of 
finding a particle in the quantum state $u_n$ on the set $\omega$. The problem of characterizing the probability measures in $\mathcal{N}(\infty)$ has attracted 
a lot of attention in the last forty years especially in the context of the so-called quantum ergodicity problem -- see e.g.~\cite{Ze09, Sa11, Non13} for 
recent surveys on that topic. Elements in $\ml{N}(\infty)$ are often called \emph{quantum limits}. In the case of $\IS^d$, it is well known that $\ml{N}(\infty)$ is contained in $\ml{N}$ which is, by definition, the closed convex hull 
(with respect to the weak-$\star$ topology) of the set of probability measures $\delta_\gamma$, where $\gamma$ is a closed geodesic of $(\IS^d,\text{Can})$. 
Recall that
$$\int_{\IS^d} a(x)\delta_\gamma(dx) = \frac{1}{2\pi}\int_0^{2\pi}a(\gamma(s))ds,$$
where the parametrization $\gamma(s)$ has unit speed. In the case where $V\equiv 0$, it was proved by Jakobson and Zelditch~\cite{JaZe99} that 
$\ml{N}(\infty)=\ml{N}$ -- the same result holds on other manifolds with positive curvature~\cite{AzM10,Ma08}. Again, it is natural to ask if this property remains true when $V$ does not identically 
vanish. This is of course related to the above observability question and we shall again answer to this question by the negative provided $V$ satisfies 
certain generic properties. In \cite{MaRi16} we showed that the answer remains negative on certain Zoll manifolds, even when $V$ vanishes.

We finally present a simple criterium relating asymptotic separation properties of the spectrum of the Schrödinger operator to the structure of the set $\ml{N}(\infty)$. We extend the proof given in \cite{MaRi16} to the case of potentials with vanishing Radon transform.

\section{Statement of the main results}

In order to state our results, we need to define the Radon transform of the potential $V$. Denote by $G(\IS^d)$ the space of closed geodesics on $\IS^d$, which is 
a smooth symplectic manifold~\cite{Be78}. Then, one can define the Radon transform of $V$ as follows:
$$\forall\gamma\in G(\IS^d),\ \ml{I}(V)(\gamma)=\int_{\IS^d} V(x)\delta_{\gamma}(dx).$$
This is a smooth function on $G(\IS^d)$ which can also be identified with a smooth $0$-homogeneous function on $T^*\IS^d-\{0\}.$ We denote by 
$\varphi_{\ml{I}(V)}^t$ the corresponding Hamiltonian flow on $T^*\IS^d-\{0\}$ which can itself be identified with an Hamiltonian flow on the symplectic 
manifold $G(\IS^d)$. We also define the second order average:
$$\ml{I}^{(2)}(V):=\ml{I}(V^2)-\frac{1}{2\pi}\int_0^{2\pi}\int_0^{t}\{V\circ\varphi^{t},V\circ\varphi^{s}\}dsdt,$$
where $\varphi^t$ denotes the geodesic flow on $S^*\IS^d$. This extends into a smooth $0$-homogeneous function on $T^*\IS^d-\{0\}$ that is invariant by the geodesic flow, and it can again be viewed as a 
function acting on $G(\IS^d)$. We denote by $\varphi_{\ml{I}^{(2)}(V)}^t$ its Hamiltonian flow.
\subsection{Observability of eigenfunctions}
 Our first result is the following
\begin{theo}\label{t:observability} Let $\omega$ be an open set in $\IS^d$. Suppose that one of the following conditions holds:

\noindent (i)  $\ml{I}(V)$ is non-constant and
$$K_{\omega,V}:=\left\{\gamma\in G(\IS^d):\forall t\in\IR,\ \varphi_{\ml{I}(V)}^t(\gamma)\cap\omega=\emptyset\right\}=\emptyset.$$
(ii) $\ml{I}(V)$ is constant and
$$K_{\omega,V}^{(2)}:=\left\{\gamma\in G(\IS^d):\forall t\in\IR,\ \varphi_{\ml{I}^{(2)}(V)}^t(\gamma)\cap\omega=\emptyset\right\}=\emptyset.$$
 Then, there exists $C_{\omega,V}>0$ such that, for every $u$ solution of~\eqref{e:eigenf}, one has
\begin{equation}\label{e:observability-eigenfunction-potential}0<C_{\omega,V}\leq\|u\|_{L^2(\omega)}^2=\int_{\omega}|u(x)|^2\operatorname{vol}(dx).\end{equation}
\end{theo}
Note that $K_{\omega,V}\subset K_{\omega}$ and that it may happen that $K_{\omega,V}=\emptyset$ while 
$K_{\overline{\omega}}$ contains a nonempty open set of closed geodesics -- see Remark~\ref{r:counterexample} below. In particular, this statement shows that observability 
for eigenfunctions may hold even if $K_{\overline{\omega}}\neq\emptyset$ provided that we choose a good $V$. This contrasts with the case 
of observability for the Schr\"odinger evolution on $\IS^d$ where $K_{\overline{\omega}}$ implies the failure of the observability property~\eqref{e:observability}. As in the classical argument of Lebeau, 
this Theorem follows from the unique continuation principle (for the case of low frequencies) and from the study of the microlocal lift of eigenfunctions (for the case 
of high frequencies). 

\begin{rema}\label{r:counterexample} Let us explain how to construct $\omega$ and $V$ such that $K_{\omega,V}=\emptyset$ while $K_{\overline{\omega}}\neq\emptyset$.
Recall first that the space of geodesics $G(\IS^2)$ can be identified with $\IS^2$~\cite[p.~54]{Be78}. This can be easily seen as follows. Take an oriented 
closed geodesic $\gamma$. It belongs to an unique $2$-plane in $\IR^3$ which can be oriented via the orientation of the geodesic, and 
$\gamma$ can be identified with the unit vector in $\IS^2$ which is directly orthogonal to this oriented $2$-plane. With that identification in mind, we also 
know from the works of Guillemin that $\ml{I}:V\in\ml{C}^{\infty}_{\text{even}}(\IS^2)\mapsto \ml{I}(V)\in\ml{C}^{\infty}_{\text{even}}(\IS^2)$ is an 
isomorphism~\cite{Gu76}. We can now explain how to construct $\omega$ and $V$. Write $\IS^2:=\left\{(x,y,z):x^2+y^2+z^2=1\right\}.$ Suppose first 
that the open set $\omega$ contains the north pole $(0,0,1)$ and that it 
does not intersect a small enough neighborhood of the equator $\Gamma=\{(x,y,0):x^2+y^2=1\}$. For instance, one can take $\omega$ to be 
equal to $\left\{(x,y,z):x^2+y^2+z^2=1\ \text{and}\ z>\eps\right\}$ with $\eps>0$ small enough. In particular, there are infinitely many geodesics 
which belongs to $K_{\overline{\omega}}\subset K_{\omega}$, i.e. the geometric control condition fails. In the space of geodesics $G(\IS^2)\simeq\IS^2$, the geodesics 
belonging to $K_{\overline{\omega}}$ correspond to a small neighborhood of the two poles $(0,0,-1)$ and $(0,0,1)$ of $\IS^2$. Hence, if one 
chooses $V$ in such a way that $\ml{I}(V)$ has no critical points in a slightly bigger neighborhood\footnote{This is possible thanks to Guillemin's result.}, 
then one finds that $K_{\omega,V}=\emptyset$. 

 
\end{rema}

\subsection{Description of $\ml{N}(\infty)$} Let us now turn to the related problem of characterizing the elements inside $\ml{N}(\infty)$. In this direction, 
we prove the following results:
\begin{theo}\label{t:quantumlimits} Let $\nu$ be a measure in $\ml{N}(\infty)$ and let  $\gamma\in G(\IS^d)$.

\noindent (a) One then has
$$d_{\gamma}\ml{I}(V)\neq 0 \Longrightarrow \nu(\gamma)=0.$$
(b) If $\ml{I}(V)$ is identically constant then:
$$d_{\gamma}\ml{I}^{(2)}(V)\neq 0 \Longrightarrow \nu(\gamma)=0.$$
In particular, whenever $\ml{I}(V)$ is non-constant or $\ml{I}(V)$ is constant but $\ml{I}^{(2)}(V)$ is not, one has
$$\ml{N}\neq\ml{N}(\infty).$$
\end{theo}
\begin{theo}\label{t:2d} If $d=2$, any $\nu$ in $\ml{N}(\infty)$ can be decomposed as follows:
 $$\nu = f \operatorname{vol} + \alpha \nu$$
where $f\in L^1(\IS^2)$, $\alpha\in[0,1]$ and $\nu$ belongs to $\ml{N}_{\operatorname{Crit}}(V)$ which is by definition the closed convex hull 
(with respect to the weak-$\star$ topology) of the set of probability measures $\delta_\gamma$, where $d_{\gamma}\ml{I}(V)=0$. 

\noindent If $\ml{I}(V)$ is constant then $\nu$ is supported on the set of critical points of $\ml{I}^{(2)}(V)$.
\end{theo}
Concerning the conclusion of Theorem \ref{t:2d}, we recall from Remark~\ref{r:counterexample} that 
$\ml{I}:V\in\ml{C}^{\infty}_{\text{even}}(\IS^2)\mapsto \ml{I}(V)\in\ml{C}^{\infty}_{\text{even}}(\IS^2)$ is an 
isomorphism. In particular, $\ml{I}(V)$ 
can always be identified with a smooth function on the real projective plane $\IR P^2$. Hence, for a generic 
choice of potential $V$, the set $\ml{N}_{\operatorname{Crit}}(V)$ is the convex hull of finitely many measures carried by closed geodesics depending only 
on $V$ -- see for instance~\cite[Sect.~3.4]{CdVVuN03} for discussions (and also examples) on critical points in this geometric framework. 

In the proofs of Theorems~\ref{t:observability}, \ref{t:quantumlimits}, and~\ref{t:2d} we will make use of classical methods from microlocal analysis 
that were originally developed for the study of the eigenvalue distribution by Duistermaat-Guillemin~\cite{DuGu75}, Weinstein~\cite{We77} and Colin de Verdi\`ere~\cite{CdV79}. 
Even if it sounds natural, it seems that the problem of characterizing 
$\ml{N}(\infty)$ in this geometric framework has not been explicitly considered in the 
literature before except when $V\equiv 0$~\cite{JaZe99, Ma08, AzM10, HPT}. We will show that these methods from microlocal analysis 
allow to obtain in a rather simple manner nontrivial results on the high frequency behaviour of Schr\"odinger eigenfunctions. 

\subsection{Relation with the study of eigenvalue distribution}

Finally, observe that the following Theorem allows to establish a relation 
between the study of $\ml{N}(\infty)$ and the level spacings:
\begin{theo}\label{t:spacing} Let $\lambda_1<\lambda_2<\lambda_3<\ldots$ be the sequence of distinct eigenvalues of $-\frac{\Delta}{2}+V$. 
Suppose that
 $$\lim_{j\rightarrow+\infty}\sqrt{\lambda_j}(\lambda_{j+1}-\lambda_j)=+\infty.$$
 Then, for every $\gamma$ in $G(\IS^d)$, $\delta_{\gamma}\in\ml{N}(\infty)$.
 
 Moreover, the same conclusion holds if $\ml{I}(V)$ is constant and if we suppose 
 $$\lim_{j\rightarrow+\infty}\lambda_j^{\frac{3}{2}}(\lambda_{j+1}-\lambda_j)=+\infty.$$
\end{theo}
The first part of this Theorem was proved in~\cite[Sect.~6]{MaRi16} in the slightly more general framework of Zoll manifolds. It follows the strategy presented in~\cite{Ma08} which 
consists in computing quantum limits using coherent states for the non-stationary Schr\"odinger equation -- see also~\cite{HPT} for a recent, different, proof of the result in~\cite{Ma08}. 
When $\ml{I}(V)$ is constant, the proof from~\cite{MaRi16} can be adapted and we shall briefly explain in paragraph~\ref{ss:spacing} which modifications should be made to get the second part. This result 
combined with Theorem~\ref{t:quantumlimits} shows that, if $\ml{I}(V)$ is non constant, then we can find a subsequence of distinct eigenvalues $(\lambda_j)_{j\in S}$ such that 
$\lambda_{j+1}-\lambda_j=\ml{O}(\lambda_j^{-\frac{1}{2}})$ for $j\in S$ tending to $+\infty$. When $\ml{I}(V)$ is constant but $\ml{I}^{(2)}(V)$ is not, then we deduce that $\lambda_{j+1}-\lambda_j=\ml{O}(\lambda_j^{-\frac{3}{2}})$. In 
other words, this gives simple criteria under which you can prove the existence of distinct eigenvalues which are asymptotically very close. In the case where $d=2$ and where $\ml{I}(V)$ is constant, 
a much stronger result on level spacings was recently
proved in~\cite{HaHiSj15}. Yet, in higher dimension or in the case of vanishing averages, it is not clear that such a 
result could be directly deduced from the classical results on the distribution of eigenvalues from~\cite{We77, Gu78a, CdV79, Ur85, Ze96}.

\subsection{The case of nonstationary solutions and general Zoll manifolds}

A natural extension of all the above problems is to consider the case of quasimodes and of nonstationary solutions when $(M,g)$ is a more general Zoll 
manifold~\cite{Be78}. These issues were discussed in great details in~\cite{MaRi16}. In order to emphasize the main geometric ideas and to avoid the technical 
issues inherent to these generalizations, we only focus here on the simpler framework described above. We refer the interested reader to 
this reference for more precise results. Note that our method combined 
with earlier results of Zelditch~\cite{Ze96, Ze97} allows in fact to show that, when $V\equiv 0$, one has $\ml{N}(\infty)\neq\ml{N}$ for many Zoll 
surfaces of revolution on $\IS^2$ -- see~\cite[Th.~1.4]{MaRi16} for the precise statement.

\section{Semiclassical measures and their invariance properties}

In order to prove the above results, we will make use of the so-called semiclassical measures~\cite{Ge91} -- see also~\cite[Chap.~5]{Zw12} for an 
introduction on that topic. In particular, we introduce the semiclassical parameter $\hbar=\lambda^{-1/2}$ and we are interested in the solutions of the 
following problem:
\begin{equation}\label{e:eigenf-semiclassical}
\left(-\frac{\hbar^2\Delta}{2}+\hbar^2V(x)\right)u_{\hbar}(x)= u_{\hbar}(x), \quad ||u_{\hbar}||_{L^2(\IS^d)}=1,
\end{equation}
in the semiclassical limit $\hbar\rightarrow 0^+$.

\subsection{Semiclassical measures} One can define the Wigner distribution of the quantum states $u_{\hbar}$:
$$\mu_{\hbar}:a\in\ml{C}^{\infty}_c(T^*\IS^d)\mapsto\la u_{\hbar},\Oph(a) u_{\hbar}\ra_{L^2(\IS^d)},$$
where $\Oph(a)$ is a pseudodifferential operator in $\Psi^{-\infty}(\IS^d)$ with principal symbol $a$ -- see~\cite[Ch.~4 and~14]{Zw12}. From the 
Calder\'on-Vaillancourt~\cite[Ch.~5]{Zw12}, the sequence $(\mu_{\hbar})_{\hbar\rightarrow 0^+}$ is bounded in 
$\ml{D}'(T^*\IS^d)$. Thus, one can extract subsequences and we denote by $\ml{M}(\infty)\subset\ml{D}'(T^*\IS^d)$ the set 
of all possible accumulation points (as $\hbar\rightarrow 0^+$) when $(u_{\hbar})_{\hbar\rightarrow 0^+}$ varies among sequences
 satisfying~\eqref{e:eigenf-semiclassical}. From the G\aa{}rding inequality~\cite[Ch.~4]{Zw12}, one can in fact verify that any 
 $\mu\in\ml{M}(\infty)$ is a finite positive measure on $T^*\IS^d$. Hence, any such $\mu$ is called a \emph{semiclassical measure}. Then, applying the composition rule  
 for pseudodifferential operators~\cite[Ch.~4]{Zw12}, one can show that any $\mu\in\ml{M}(\infty)$ is supported on the unit cotangent bundle $S^*\IS^d$
and that
\begin{equation}\label{e:pushforward}
 \ml{N}(\infty):=\left\{\int_{S^*_x\IS^d}\mu(x,d\xi):\ \mu\in\ml{M}(\infty)\right\}.
\end{equation}
For more details on these facts, we refer the reader to~\cite[Ch.~5]{Zw12}. Finally, as a warm up, let us briefly remind how to prove that 
these measures are invariant by 
the geodesic flow $\varphi^s$, i.e. the Hamiltonian flow associated with the function $\frac{\|\xi\|^2_x}{2}$. For that purpose, we write, given 
$u_{\hbar}$ satisfying~\eqref{e:eigenf-semiclassical} and any $a$ in $\ml{C}^{\infty}_c(T^*\IS^d)$,
\begin{equation}\label{e:eigenvalue-equation}
\left\la u_{\hbar},\left[-\frac{\hbar^2\Delta}{2}+\hbar^2V,\Oph(a)\right]u_{\hbar}\right\ra=0.
\end{equation}
We can now apply the commutation rule for pseudodifferential operators~\cite[Ch.~4]{Zw12} combined with the Calder\'on-Vaillancourt Theorem:
$$\frac{\hbar}{i}\left\la u_{\hbar},\Oph\left(\left\{\frac{\|\xi\|^2}{2},a\right\}\right)u_{\hbar}\right\ra=\ml{O}(\hbar^2),$$
where $\{,\}$ is the Poisson bracket. Dividing this equality by $\hbar$ and letting $\hbar$ goes to $0$ in this equality, one finds (after a 
possible extraction) that $\mu(\{\|\xi\|^2,a\})=0$ for every $a$ in $\ml{C}^{\infty}_c(T^*\IS^d)$. From the properties of $\mu$, it exactly shows that 
elements in $\ml{M}(\infty)$ are invariant by the geodesic flow $\varphi^s$ acting on $S^*\IS^d$.

\subsection{Weinstein's averaging method} Note that all the arguments so far are valid on a general compact Riemannian manifold and we shall now see which 
extra properties can be derived in the case of $\IS^d$ endowed with its canonical metric. For that purpose, we need to fix some conventions and to 
collect some well-known facts on the spectral properties of the Laplace-Beltrami operator on $\IS^d$. First, given any 
$a$ in $\ml{C}^{\infty}_c(T^*\IS^d-\{0\})$, 
we introduce the Radon transform of $a$:
$$\ml{I}(a)(x,\xi):=\frac{1}{2\pi}\int_0^{2\pi}a\circ\varphi^{s\|\xi\|_{x}}(x,\xi)ds.$$
In the case of $V$, this definition can be identified with the Radon transform that were defined in the introduction. 
We will now define the equivalent of this operator at the quantum level following the seminal work of 
Weinstein~\cite{We77} -- see also~\cite{DuGu75, CdV79} for more general geometric frameworks. 
Recall that the eigenvalues of $-\Delta$ are of the form 
$$\lambda_k=\left(k+\frac{d-1}{2}\right)^2-\frac{(d-1)^2}{4},$$
where $k$ runs over the set of nonnegative integer. In particular, we can write
\begin{equation}\label{e:decompose-Laplace}
 -\Delta= A^2-\left(\frac{d-1}{2}\right)^2,
\end{equation}
where $A$ is a selfadjoint pseudodifferential operator of order $1$ with principal symbol $\|\xi\|_x$ and satisfying
\begin{equation}\label{e:period-quantum}
e^{2i\pi A}=e^{i\pi(d-1)}\text{Id}.
\end{equation}
Given $a$ in $\ml{C}^{\infty}_c(T^*\IS^d-\{0\})$, we then set, by analogy with the Radon transfom of $a$,
$$\ml{I}_{\text{qu}}(\Oph(a)):=\frac{1}{2\pi}\int_0^{2\pi}e^{-is A}\Oph(a)e^{isA}ds.$$
An important observation which seems to be due to Weinstein~\cite{We77} is that the following exact commutation relation holds:
$$\left[\ml{I}_{\text{qu}}(\Oph(a)),A\right]=0.$$
In particular, from~\eqref{e:decompose-Laplace}, one has
\begin{equation}\label{e:commute}
\left[\ml{I}_{\text{qu}}(\Oph(a)),\Delta\right]=0.
\end{equation}
Finally, the Egorov Theorem allows to relate the operator $\ml{I}_{\text{qu}}(\Oph(a))$ to the classical Radon transform as follows:
\begin{equation}\label{e:egorov}
 \ml{I}_{\text{qu}}(\Oph(a))=\Oph(\ml{I}(a))+\hbar R,
\end{equation}
where $R$ is a pseudodifferential operator in $\Psi^{-\infty}(\IS^d)$

\subsection{Extra invariance properties on $\IS^d$} Let us now apply these properties to derive some invariance properties of the elements in 
$\ml{M}(\infty)$. We fix $\mu$ in $\ml{M}(\infty)$ which is generated by a sequence $(u_{\hbar})_{\hbar\rightarrow 0^+}$ and 
$a$ in $\ml{C}^{\infty}_c(T^*\IS^d-\{0\})$. We rewrite~\eqref{e:eigenvalue-equation} with $\ml{I}_{\text{qu}}(\Oph(a))$ 
instead of $\Oph(a)$. According to~\eqref{e:commute}, this implies that
$$\left\la u_{\hbar},\left[V,\ml{I}_{\text{qu}}(\Oph(a))\right]u_{\hbar}\right\ra=0.$$
Combining~\eqref{e:egorov} with the commutation formula for pseudodifferential operators and the Calder\'on-Vaillancourt theorem, we find then
$$\frac{\hbar}{i}\left\la u_{\hbar},\Oph\left(\{V,\ml{I}(a)\}\right)u_{\hbar}\right\ra=\ml{O}(\hbar^2).$$
Hence, after letting $\hbar$ goes to $0$, one finds that
$$\mu(\{V,\ml{I}(a)\})=0.$$
Applying the invariance by the geodesic flow, one finally gets that 
\begin{equation}\label{e:main-result}\mu(\{\ml{I}(V),a\})=0.\end{equation}
This is valid for any smooth test function $a$ in $\ml{C}^{\infty}_c(T^*\IS^d-\{0\})$. Thus, we have just proved that \textbf{any $\mu$ in $\ml{M}(\infty)$ 
is invariant by the Hamiltonian flow $\varphi_{\ml{I}(V)}^t$ of $\ml{I}(V)$} which is well defined on $S^*\IS^d\subset T^*\IS^d-\{0\}$. In other words, any element in 
$\ml{M}(\infty)$ is an invariant measure for the system 
$$F:T^*\IS^{d}-\{0\}\ni (x,\xi)\mapsto \left(\frac{\|\xi\|_x^2}{2},\ml{I}(V)(x,\xi)\right)\in\IR^2,$$
which is completely integrable in dimension $2$. Theorems~\ref{t:quantumlimits} and \ref{t:2d} follow then from classical arguments on 
integrable systems -- see e.g. paragraph~3.3 in~\cite{MaRi16} for part (a) of Theorem \ref{t:quantumlimits} and Corollary~4.4 of that reference for the first conclusion of Theorem \ref{t:2d}.

\subsection{The case of vanishing averages}\label{ss:vanishing} One can easily observe that the results we have proved so far are empty if we suppose that $V$ is an odd 
function on $\IS^d$. This is due to the fact that identity~\eqref{e:main-result} does not provide any non-trivial information on $\mu$ in that case. We would now like to explain how one can obtain a new invariance relation in 
that case -- namely invariance by the Hamiltonian flow of the second order average $\ml{I}^{(2)}_V$. This is enough to prove part (b) of Theorem \ref{t:quantumlimits} and complete the proof of Theorem \ref{t:2d}. This problem was not considered in~\cite{MaRi16} and we will briefly 
expose how some ideas of Guillemin and Uribe~\cite{Gu78a, Ur85} can be applied to treat this case. 

From this point on, we suppose that $V$ is \emph{an 
odd function} on $\IS^d$. In particular, its Radon 
transform identically vanishes. Recall also from~\cite[Lemma~3.1]{Gu78a} that its quantum counterpart also identically vanishes, i.e.
\begin{equation}\label{e:guillemin}\ml{I}_{\text{qu}}(V):=\frac{1}{2\pi}\int_0^{2\pi}e^{-isA}Ve^{isA}ds=0.\end{equation}
whenever $V$ is an odd function. Following~\cite{Ur85} and given a bounded operator $C$ on $L^2(\IS^d)$, one can define 
$$\ml{I}_{\text{qu}}(C):=\frac{1}{2\pi}\int_0^{2\pi}e^{-isA}Ce^{isA}ds,$$
and 
$$\sigma(C):=-\frac{1}{2\pi}\int_0^{2\pi}dt\left(\int_0^te^{-isA}Ce^{isA}ds\right).$$
As was already observed, one has $[A,\ml{I}_{\text{qu}}(C)]=0$. For $\sigma(C)$, the following holds:
\begin{equation}\label{e:sigma-commute}
 [A,\sigma(C)]=i(C-\ml{I}_{\text{qu}}(C)).
\end{equation}
We now set
$$U_{\hbar}(t):=\exp\left(-it\hbar\sigma(Q)\right),$$
where $Q_\hbar$ is an $\hbar$-pseudodifferential operator in $\Psi^{-1}(\IS^d)$ that has to be determined. We also fix $a$ in $\ml{C}^{\infty}_c(T^*\IS^d-\{0\})$.
In order to motivate the upcoming calculation, we now write~\eqref{e:eigenvalue-equation} with 
$U_{\hbar}(-1)\ml{I}_{\text{qu}}(\Oph(a))U_{\hbar}(1)$ instead of $\Oph(a)$:
$$\left\la u_{\hbar},\left[\frac{\hbar^2A^2}{2}+\hbar^2V,U_{\hbar}(-1)\ml{I}_{\text{qu}}(\Oph(a))U_{\hbar}(1)\right]u_{\hbar}\right\ra=0.$$
Equivalently, this can be rewritten as
\begin{equation}\label{e:conjugated-test-operator}\left\la u_{\hbar},U_{\hbar}(-1)\left[U_{\hbar}(1)\left(\frac{\hbar^2A^2}{2}+\hbar^2V\right)U_{\hbar}(-1),\ml{I}_{\text{qu}}(\Oph(a))\right]U_{\hbar}(1)u_{\hbar}\right\ra=0.\end{equation}
Using the fact that $V$ is odd, we would now like to choose an appropriate $Q_\hbar$ such that
$$U_{\hbar}(1)\left(\frac{\hbar^2A^2}{2}+\hbar^2V\right)U_{\hbar}(-1)=\frac{\hbar^2A^2}{2}+\hbar^4 Q_\hbar^1,$$
for some bounded pseudodifferential operator $Q_\hbar^1$ to be determined. This kind of normal form for the Schr\"odinger operator on $\IS^d$ was for 
instance obtained by Guillemin in~\cite[Sect.~3]{Gu78a} and by Uribe in~\cite[Sect.~4 and 6]{Ur85}. Let us recall their argument.

We first use~\eqref{e:sigma-commute} and the composition formula for pseudodifferential operators 
to write that
$$U_{\hbar}(1)\hbar AU_{\hbar}(-1)=\hbar A-\hbar^2(Q_\hbar-\ml{I}_{\text{qu}}(Q_\hbar))+i\frac{\hbar^3}{2}[\sigma(Q_\hbar),Q_\hbar-\ml{I}_{\text{qu}}(Q_\hbar)]+\ml{O}_{\Psi^{-1}(\IS^d)}(\hbar^5).$$
Hence, if we square this expression, we find, using the composition formula one more time,
\begin{eqnarray*}
U_{\hbar}(1)\hbar^2 A^2U_{\hbar}(-1) & = &\hbar^2A^2-\hbar^2 [\hbar A(Q_\hbar-\ml{I}_{\text{qu}}(Q_\hbar))-(Q_\hbar-\ml{I}_{\text{qu}}(Q_\hbar))\hbar A]\\
 & &+i\frac{\hbar^3}{2} \left(\hbar A[\sigma(Q_\hbar),Q_\hbar-\ml{I}_{\text{qu}}(Q_\hbar)]+[\sigma(Q_\hbar),Q_\hbar-\ml{I}_{\text{qu}}(Q_\hbar)]\hbar A\right)\\
 & & +\hbar^4(Q_\hbar-\ml{I}_{\text{qu}}(Q_\hbar))^2+\ml{O}_{\Psi^{0}(\IS^d)}(\hbar^5).
\end{eqnarray*}
Similarly, one has
$$U_{\hbar}(1)\hbar^2 VU_{\hbar}(-1)=\hbar^2V-i\hbar^3[\sigma(Q_\hbar),V]+\ml{O}_{\Psi^{0}(\IS^d)}(\hbar^5).$$
Thus, if we want to cancel the term $\hbar^2 V$ in~\eqref{e:conjugated-test-operator}, we have to impose that $(\hbar A)(Q_\hbar-\ml{I}_{\text{qu}}(Q_\hbar))+
(Q_\hbar-\ml{I}_{\text{qu}}(Q_\hbar))(\hbar A)$ is equal to $2 V$ (at least at first order). Define, for every  bounded operator $B$, 
$$\ml{L}(B)=B A^{-1}+A^{-1}B,$$
and set finally
$$Q_\hbar=\frac{1}{\hbar}\ml{L}\left(2V- \frac{A V A^{-1}+ A^{-1} V A}{2}\right),$$
which is in $\Psi^{-1}(\IS^d)$ with a principal symbol equal to $q(x,\xi)=V(x)/\|\xi\|_x$: 
$$Q_\hbar=\Oph(q)+\hbar R_\hbar,\quad \text{ with } R_\hbar \text{ bounded.}$$
Observe that, as $V$ is odd, one can verify that $\ml{I}(q)\equiv 0$. In the 
following, we will denote by $\sigma(q)$ the principal symbol of the operator $\sigma(Q_\hbar)$. Thanks to~\eqref{e:guillemin}, we also know that $\ml{I}_{\text{qu}}(Q_\hbar)=0$ from which we can infer
$$(\hbar A)(Q_\hbar-\ml{I}_{\text{qu}}(Q_\hbar))+(Q_\hbar-\ml{I}_{\text{qu}}(Q_\hbar))(\hbar A)=(\hbar A)Q_\hbar+Q_\hbar(\hbar A).$$
In other words, it remains to compute the difference between $(\hbar A)Q_\hbar+Q_\hbar(\hbar A)$ and $2V$:
$$(\hbar A)Q_\hbar+Q_\hbar(\hbar A)-2V=A^{-1}VA+AVA^{-1}-V-\frac{A^{-2}VA^2+A^2VA^{-2}}{2}.$$
We now write that $ABA^{-1}=[A,BA^{-1}]+B$ and $A^{-1}BA=[A^{-1}B,A]+B$ which implies
$$(\hbar A)Q_\hbar+Q_\hbar(\hbar A)-2V=-\frac{[A,[A,VA^{-1}]A^{-1}]+[A^{-1}[A^{-1}V,A],A]}{2}.$$
According to the composition rules for pseudodifferential operators, we 
find that $(\hbar A)Q_\hbar+Q_\hbar(\hbar A)-2V$ belongs to $\hbar^2\Psi^0(\IS^d)$ with principal symbol equal to
$$\hbar^2r(x,\xi)=-\hbar^2\frac{\{\|\xi\|_x,\{\|\xi\|_x,V(x)\|\xi\|^{-1}_x\}\|\xi\|^{-1}_x\}+
\{\|\xi\|_x^{-1}\{V(x)\|\xi\|^{-1}_x,\|\xi\|_x\},\|\xi\|_x\}}{2}.$$
Note that, as $V$ is odd, one can verify that $\ml{I}(r)\equiv 0$. Combining these equalities, we find that
$$U_{\hbar}(1)\left(\frac{\hbar^2A^2}{2}+\hbar^2V\right)U_{\hbar}(-1)=\frac{\hbar^2A^2}{2}+\hbar^4 \Oph\left(q_1\right)+\ml{O}_{\Psi^{0}(\IS^d)}(\hbar^5),$$
where 
$$q_1(x,\xi):=\frac{q(x,\xi)^2+\|\xi\|_x\{\sigma(q),q\}(x,\xi)-2\{\sigma(q),V\}(x,\xi)-r(x,\xi)}{2}.$$
Insert this identity in~\eqref{e:conjugated-test-operator} and apply~\eqref{e:commute} and ~\eqref{e:egorov} to derive that
\begin{equation}\label{e:conjugated-test-operator-2}\hbar^5\left\la u_{\hbar},\Oph(\{q_1,\ml{I}(a)\})u_{\hbar}\right\ra=\ml{O}(\hbar^6).\end{equation}
If we let $\hbar$ go to $0$, we find that the corresponding semiclassical measure $\mu$ verifies $\mu(\{q_1,\ml{I}(a)\})=0$. From the invariance of 
$\mu$ by the geodesic flow and from the relation $\ml{I}(r)\equiv 0$, this implies that
\begin{equation}\label{e:invariance-odd}
 \mu(\{q(x,\xi)^2+\|\xi\|_x\{\sigma(q),q\}(x,\xi)-2\{\sigma(q),V\}(x,\xi),\ml{I}(a)\})=0.
\end{equation}
Then, use that $\{\|\xi\|_x,\ml{I}(a)\}=0$ and that $\mu$ is supported in $S^*\IS^d$ in order to show that this is equivalent to
$$\mu(\{V^2-\{\sigma(V),V\},\ml{I}(a)\})=0.$$
Using the invariance of $\mu$ by the geodesic flow, we find that
\begin{equation}\label{e:invariance-odd-end}
 \mu\left(\left\{\ml{I}(V^2)-\frac{1}{2\pi}\int_0^{2\pi}\int_0^t\{V\circ\varphi^t,V\circ\varphi^s\}dsdt,a\right\}\right)=0,
\end{equation}
which replaces~\eqref{e:main-result} when $V$ is an odd function. 

\subsection{Observability estimates} Let us now give the proof of Theorem~\ref{t:observability}. Suppose by contradiction that this result is not true. It means 
that there exists a sequence $(u_n)_{n\geq 1}$ of solutions of~\eqref{e:eigenf} such that $\|u_n\|_{L^2(\omega)}\rightarrow 0$. From the unique continuation 
principle (see e.g.~\cite{LRLe12}) and using the fact that $\omega$ is a non empty open set, one can verify that $\lambda_n$ has to converge to infinity. Up to an extraction, we can 
suppose that $(u_n)_{n\geq 1}$ generates an unique semiclassical measure $\mu$. Using the invariance by the geodesic flow, one knows that 
$\mu(S^*\omega)=\mu(\ml{I}(\mathbf{1}_{\omega}))$. Suppose now that $\ml{I}(V)$ is non-constant. Then, as $\mu$ is a postive measure, one knows that
$$\lim_{n\rightarrow+\infty}\|u_n\|_{L^2(\omega)}^2\geq\mu(S^*\omega)\geq \frac{1}{T}\inf_{\rho\in S^*\IS^d}\int_0^T\ml{I}(\mathbf{1}_{\omega})\circ\varphi_{\ml{I}(V)}^s(\rho)ds.$$
From the fact that $K_{\omega,V}=\emptyset$, one knows that, for $T>0$ large enough, this lower bound is positive which implies the expected contradiction as 
the upper bound vanishes by hypothesis. When $\ml{I}(V)$ is constant it suffices to reproduce this argument using the invariance of semiclassical measures by the Hamiltonian flow of $\ml{I}^{(2)}(V)$.

\subsection{Relation to eigenvalue distribution}\label{ss:spacing}

In this last paragraph, we briefly explain the main lines of the proof of Theorem~\ref{t:spacing}. We only treat the second part of the Theorem which was not discussed in~\cite{MaRi16}. Therefore, as in paragraph~\ref{ss:vanishing}, we suppose that $V$ is odd. 

First of all, fix a point $(x_0,\xi_0)$ in $S^*\IS^d$ and a normalized sequence 
$(u_{\hbar}^{x_0,\xi_0})_{\hbar\rightarrow 0^+}$ of coherent states whose semiclassical measure is $\delta_{x_0,\xi_0}$. Recall 
from~\cite[Sect.~6.1]{MaRi16} that, up to some spectral truncation, we can always suppose that 
\begin{equation}\label{e:spectral-decomp}
u_{\hbar}^{x_0,\xi_0}=\sum_{\{j:\frac{1}{4}\leq \lambda_j\hbar^2\leq 1\}}c_{\hbar}^{x_0,\xi_0}(j)\hat{v}_{\hbar}^{x_0,\xi}(j), 
\end{equation}
where, for every choice of parameters, $c_{\hbar}^{x_0,\xi_0}(j)\geq 0$ and $\hat{v}_{\hbar}^{x_0,\xi_0}(j)$ is normalized in $L^2(\IS^d)$ and verifies
$$\left(-\frac{\Delta}{2}+V\right)\hat{v}_{\hbar}^{x_0,\xi_0}(j)=\lambda_j\hat{v}_{\hbar}^{x_0,\xi_0}(j).$$ 
We now let
$(\tau_{\hbar})_{\hbar\rightarrow 0^+}$ be a sequence of times such that
\begin{equation}\label{e:level-spacing}\lim_{\hbar\rightarrow 0^+}\tau_{\hbar}\min\left\{\lambda_{j+1}-\lambda_j:\ \frac{1}{4}\leq \hbar^2\lambda_j\leq 1\right\}=+\infty.\end{equation}
From our assumption on the level spacing, we can in fact suppose that $\tau_{\hbar}=o(\hbar^{-3})$.
As in the previous sections, we fix $a$ in $\ml{C}^{\infty}_c(T^*\IS^d-\{0\})$ and we consider the time dependent Wigner distribution:
$$\la \mu_{\hbar}^{x_0,\xi_0}(t),a\ra:=\la v_{\hbar}^{x_0,\xi_0}(t\tau_{\hbar}),U_{\hbar}(-1)\ml{I}_{\text{qu}}(\Oph(a))U_{\hbar}(1)
v_{\hbar}^{x_0,\xi_0}(t\tau_{\hbar})\ra,$$
where $v_{\hbar}^{x_0,\xi_0}(t\tau_{\hbar})$ is the solution at time $t\tau_{\hbar}$ of~\eqref{e:nonsemiclassical-schr} with 
initial condition $u_{\hbar}^{x_0,\xi_0}$. If we differentiate this expression with respect to time and if we argue as in paragraph~\ref{ss:vanishing}, we find 
that
$$\frac{d}{dt}\la \mu_{\hbar}^{x_0,\xi_0}(t),a\ra= \ml{O}(\hbar^3\tau_{\hbar}).$$
Recall that, for a general $V$, the proof from~\cite{MaRi16} gave a remainder term of order $\ml{O}(\tau_{\hbar}\hbar)$ and that we did not  
introduce the Fourier integral operator $U_{\hbar}(1)$ in our argument. Integrating this expression 
between $0$ and $t$ and using our assumption that $\tau_{\hbar}=o(\hbar^{-3})$, we find
$$\la \mu_{\hbar}^{x_0,\xi_0}(t),a\ra=\ml{I}(a)(x_0,\xi_0)+o(1).$$
If we now fix $\theta$ in $\ml{S}(\IR)$ whose Fourier transform $\ml{F}(\theta)$ is compactly 
supported and verifies $\ml{F}(\theta)(0)=1$, then we find that
$$\int_{\IR}\theta(t)\la\mu_{\hbar}^{x_0,\xi_0}(t),a\ra dt=\ml{I}(a)(x_0,\xi_0)+o(1).$$
Using the spectral decomposition~\eqref{e:spectral-decomp} and~\eqref{e:level-spacing}, we obtain the following averaging formula:
$$\sum_{\{j:\frac{1}{4}\leq \lambda_j\hbar^2\leq 1\}}c_{\hbar}^{x_0,\xi_0}(j)^2\la \hat{v}_{\hbar}^{x_0,\xi_0}(j),U_{\hbar}(-1)\ml{I}_{\text{qu}}(\Oph(a))U_{\hbar}(1) 
\hat{v}_{\hbar}^{x_0,\xi_0}(j)\ra=\ml{I}(a)(x_0,\xi_0)+o(1),$$
which yields after simplification
$$\sum_{\{j:\frac{1}{4}\leq \lambda_j\hbar^2\leq 1\}}c_{\hbar}^{x_0,\xi_0}(j)^2\la \hat{v}_{\hbar}^{x_0,\xi_0}(j),\Oph(a) 
\hat{v}_{\hbar}^{x_0,\xi_0}(j)\ra=\ml{I}(a)(x_0,\xi_0)+o(1),$$
Recall that, as $u_{\hbar}^{x_0,\xi_0}$ was chosen to be normalized in $L^2(\IS^d)$, one has $\sum_{j}c_{\hbar}^{x_0,\xi_0}(j)^2=1$. Arguing as in the proof of the Quantum Ergodicity Theorem -- see~\cite[Sect.~6.3]{MaRi16} for details, we can obtain the following variance estimate:
\begin{equation}\label{e:variance}
 \sum_{\{j:\frac{1}{4}\leq \lambda_j\hbar^2\leq 1\}}c_{\hbar}^{x_0,\xi_0}(j)^2\left|\left\la \hat{v}_{\hbar}^{x_0,\xi_0}(j),\Oph(a) \hat{v}_{\hbar}^{x_0,\xi_0}(j)\right\ra-\ml{I}(a)(x_0,\xi_0)\right|^2=o(1),
\end{equation}
which is sufficient to conclude the proof of the Theorem thanks to the 
Bienaym\'e-Tchebychev Theorem -- see~\cite[Sect.~6.4]{MaRi16} for details.

\section*{Acknowledgements}

The present note has been written for the proceedings of the workshop \emph{Probabilistic Methods in Spectral Geometry and PDE} which were held in the CRM of 
Montr\'eal at the end of the summer 2016. The authors thank the organizers of this meeting for the opportunity to expose their work~\cite{MaRi16} and some 
further developments of it in these proceedings.

\end{document}